\documentstyle[12pt]{article}
\newcommand{\const}{\mathop{\rm const}\limits}

\newcommand{\mod}{\mathop{\rm mod}\limits}

\newcommand{\Var}{\mathop{\rm Var}\limits}

\begin{document}

\begin{center}

{\bf SHARP MOMENT ESTIMATES FOR }\\

\vspace{3mm}

{\bf POLYNOMIAL MARTINGALES }\\

\vspace{4mm}

{\bf E. Ostrovsky, L.Sirota.}\\

\vspace{5mm}

Department of Mathematics, Bar-Ilan University, Ramat-Gan, 59200,Israel.\\
e-mail: \ eugostrovsky@list.ru \\

\vspace{3mm}

Department of Mathematics, Bar-Ilan University, Ramat-Gan, 59200,Israel.\\
e-mail: \ sirota3@bezeqint.net \\

\vspace{4mm}

{\bf Abstract.} \par

\vspace{3mm}

\end{center}

  In this paper non-asymptotic  moment
estimates are derived for  tail of distribution for discrete time polynomial  martingale
by means of martingale differences as a rule in the terms of {\it unconditional} and
{\it unconditional relative} moments and tails of distributions of summands. \par
 We show also the exactness of obtained estimations. \par

\vspace{4mm}

 \par

\vspace{3mm}


\vspace{4mm}

\section{ Introduction. Notations. Statement of problem. Announce.} \par

\vspace{3mm}

 Let $ (\Omega,F,{\bf P} ) $ be a probabilistic space, which will be presumed sufficiently
rich when we construct  examples (counterexamples),
$ \xi(i,1), \xi(i,2), \ldots,\xi(i,n), \ n \le \infty $ being a {\it family} of a centered
$ ({\bf E} \xi(i,m) = 0, i=1,2,\ldots,n) $ martingale-differences on the basis of the fixed
{\it flow }of $ \sigma - $ fields (filtration) $ F(i): F(0) =
\{\emptyset, \Omega \}, \ F(i) \subset F(i+1) \subset F, \
\xi(i,0) := 0; \ {\bf E} |\xi(i,m)| < \infty, $ and for every $ i \ge 0,
 \ \forall k = 0,1,\ldots,i -1 \ \Rightarrow  $

$$
{\bf E} \xi(i,m)/F(k) = 0; \ {\bf E}\xi(i,m)/F(i) = \xi(i,m) \  (\mod \ {\bf P}). \eqno(1.0)
$$

\vspace{3mm}

 Let also $ I = I(n) = I(d; n) = \{i_1; i_2; \ldots; i_d \} $ be the set of indices of the
form $ I(n) = I(d; n) = \{ \vec{i} \}  = \{i \} = \{i_1, i_2, \ldots, i_d \} $ such that $ 1 \le i_1 < i_2 <
i_3 < i_{d-1} < i_d \le n; \  J = J(n) = J(d; n) $ be the set of indices of
the form (subset of  $ I(d; n))  \  J(d; n) = J(n) =   \{ \vec{ j} \} = \{ j  \} = \{j_1; j_2; \ldots; j_{d-1} \} $
such that $ 1 \le j_1 < j_2 \ldots < j_{d-1} \le n- 1; \ b(i) = \{ b(i_1; i_2; \ldots i_d) \} $ be a
$ d  $ dimensional numerical non-random sequence symmetrical relative to all argument permutations,

$$
\vec{i} \in I \  \Rightarrow \xi(\vec{i}) \stackrel{def}{=} \prod_{s=1}^d \xi(i_s, s);\eqno(1.1a)
$$

$$
\vec{j} \in J \  \Rightarrow \xi(\vec{j}) \stackrel{def}{=} \prod_{s=1}^{d-1} \xi(i_s, s); \eqno(1.1b)
$$

$$
\sigma^2(i_k,s) := \Var(\xi_{i_k,s}),
$$

$$
\sigma^2(\vec{i}) := \prod_{s=1}^d \sigma^2 ({i_s,s}), \ \vec{i} = \{i_1, i_2, \ldots, i_d \} \in I;
$$

$$
\sigma^2(n,\vec{j}):= \prod_{s=1}^{d-1} \sigma^2 ({j_s,s}), \ \vec{j} = \{j_1, j_2, \ldots, j_{d-1} \} \in J;
$$

$$
Q_d = Q(d,n, \{  \xi(\cdot) \}  ) = Q(d,n) = Q(d,n, \vec{b}) = \sum_{\vec{i} \in I(d,n) } b(\vec{i}) \xi(\vec{i}) \eqno(1.2)
$$
being a homogeneous polynomial (random polynomial) of power $ d $ on the
random variables $ \xi(\cdot, \cdot) $  "without diagonal members", (on the other
words, multiply stochastic integral over discrete stochastic martingale
measure), $ n $ be an integer number: $ n = 1, 2, \ldots, $ in the case $ n = \infty $
we will understood $ Q(d, \infty) $  as a limit $ Q(d;\infty) = \lim_{n \to \infty} Q(d; n), $ if
there exists with probability one.\par
 Here $ b = \vec{b} = b(i), \ i \in I(d,n) $ be arbitrary non-random numerical sequence, and we denote

$$
||b||^2 = \sum_{i \in I(d,n)} b^2(i); \ b \in B \hspace{6mm} \Leftrightarrow ||b|| = 1. \eqno(1.2b)
$$

 We will denote also in the simple case when $ \vec{b} = \vec{b}_d = \vec{1} = \{ 1,1,\ldots,1  \} \ d \ - $ times

$$
R(d) =  Q(d,n, \vec{1}) = \sum_{\vec{i} \in I(d,n) }  \xi(\vec{i}). \eqno(1.3)
$$

 It is obvious that the sequence $ (Q(d; n); F(n)); n = 1, 2, 3, \ldots, $ is a
martingale ("polynomial martingale"). \par

The case of non-homogeneous polynomial is considered analogously. \par

We denote as usually the $ L(p) $ norm of the r.v. $ \eta $ as follows:

 $$
 |\eta|_p = \left[ {\bf E} |\eta|^p  \right]^{1/p}, \ p > 2;
 $$
the case $  p = 2 $ is trivial for us.\par

\vspace{4mm}

{\bf We will derive the moment estimations of a Khintchine  form }

$$
U(p;d,n) = U(p) \stackrel{def}{=}
\sup_{b \in B} | Q(d,n, b, \{ \xi(\cdot) \}) |_p \le \overline{Q}(p;d,n) = \overline{Q}(p)
$$
{\bf for martingale (and following in the independent case) in the terms of unconditional moments, more
exactly, in the $ L(p) $ norms of summands:}

$$
\mu_m(p) \stackrel{def}{=} \sup_i |\xi(i,m)|_p. \eqno(1.4)
$$

 Denote also in the martingale case

 $$
 V(p) = V_d(p) \stackrel{def}{=} \prod_{m=1}^d \mu_m(d \cdot p). \eqno(1.5a)
 $$
and  for the independent variables $ \{  \xi(i,m) \} $

 $$
 W(p) = W_d(p) \stackrel{def}{=} \prod_{m=1}^d \mu_m(p). \eqno(1.5b)
 $$
 Note that if all the functions $ p \to \mu_m(p), \ m =1,2,\ldots,d $ are regular varying:

$$
\sup_{p \ge 2} [ \mu_m(d \cdot p) /\mu_m(p)] < \infty,
$$
then $ V(p) \asymp W(p). $ \par

\vspace{3mm}

 As we knew, the previous result in this direction is obtained in the article \cite{Ostrovsky4}:

$$
U(p;d,n) \le C_1(d) \cdot p^d \ V(p),
$$
 and in the independent case, i.e. when all the (two-dimensional indexed) centered  r.v. $ \xi(i,m) $
are common independent,

$$
U(p;d,n) \le C_2(d) \cdot (p^d/\ln p) \cdot W(p), \ p \ge 2.
$$

 {\it We intend to improve both these estimates to the following non-improvable as } $ p \to \infty: $

$$
U(p;d,n)  \le \gamma(d) \cdot \frac{p^d}{(\ln p)^d } \cdot V_d(p), \eqno(1.6a)
$$
in general (martingale) case  and

$$
U(p;d,n) \le \kappa(d) \cdot \frac{p^d}{(\ln p)^d } \cdot W_d(p), \eqno(1.6b)
$$
for the independent variables; it will be presumed obviously the finiteness of  $ V(p) $ and $  W(p). $ \par

 \vspace{3mm}

 There are many works about this problem; the  next list is far from being complete:
\cite{Burkholder1}, \cite{Burkholder2},
 \cite{Fan1}, \cite{Grama1}, \cite{Grama2}, \cite{Hitczenko1}, \cite{Lesign1}, \cite{Liu1}, \cite{Osekowski1}, \cite{Ostrovsky0},
\cite{Ostrovsky7}, \cite{Ostrovsky8},  \cite{Peshkir1}, \cite{Ratchkauskas1}, \cite{Volny1} etc. \par
  See also the reference therein.\par

 Notice that in the articles \cite{Lesign1}, \cite{Liu1}, \cite{Volny1}  and in many others are described some new
applications of these estimates: in the theory of dynamical system,  in the theory of polymers etc.\par

 \vspace{4mm}

\section{ Main result: moment estimation for polynomial martingales. }

\vspace{4mm}

 We must describe some new notations.  The following function was introduced by A.Osekowski (up to factor 2)
 in the article \cite{Osekowski1}:

$$
Os(p) \stackrel{def}{=} 4 \  \sqrt{2} \cdot \left( \frac{p}{4} + 1 \right)^{1/p} \cdot \left( 1 + \frac{p}{\ln (p/2)}  \right). \eqno(2.1)
$$
 Note that

$$
K = K_{Os} \stackrel{def}{=} \sup_{p \ge 4} \left[\frac{Os(p)}{p/\ln p} \right] \approx 15.7858, \eqno(2.2)
$$
the so-called Osekowski's constant. \par
 Let us define the following numerical sequence $ \gamma(d), \ d = 1,2,\ldots:  \ \gamma(1) := K_{Os} = K,   $ (initial condition) and
by the following recursion

$$
\gamma(d+1) = \gamma(d) \cdot K_{Os} \cdot \left( 1 + \frac{1}{d} \right)^d. \eqno(2.3)
$$
 Since

$$
\left( 1 + \frac{1}{d} \right)^d \le e,
$$
we conclude

$$
\gamma(d) \le K_{Os}^d \cdot e^{d-1}, \ d = 1,2,\ldots. \eqno(2.4)
$$

\vspace{3mm}

{\bf Theorem 2.1.} {\it  Let the sequence   } $ \gamma(d) $ {\it  be defined in (2.3). Then  }

$$
U(p;d,n) \le \gamma(d) \cdot \frac{p^d}{(\ln p)^d } \cdot V_d(p) =
\gamma(d) \cdot \frac{p^d}{(\ln p)^d } \cdot \prod_{m=1}^d \mu_m(d \ p). \eqno(2.5)
$$

\vspace{3mm}

{\bf Proof.} \\

\vspace{3mm}

{\bf 0.} We will use the induction method over the "dimension" $  d, $ as in the article of authors \cite{Ostrovsky4},
starting from the value $ d=1. $ \\

{\bf 1.} One dimensional case  $  d=1.$ We apply the celebrate result belonging to A.Osekowski \cite{Osekowski1}:

$$
\left|\sum_{k=1}^n \xi_k  \right|_p \le C_{Os}(p) \cdot
\left\{ \left| \left( \sum_{k=1}^n {\bf E} \xi_k^2/F(k-1) \right)^{1/2}    \right|_p +
\left| \left( \sum_{k=1}^n |\xi_k|^p  \right)^{1/p}   \right|_p  \right\} \stackrel{def}{=} \eqno(2.6)
$$

$$
C_{Os}(p)  \left\{  S_1(p) + S_2(p) \right\}, \eqno(2.6a)
$$
in our notations; $  \xi_k = \xi(1,k), \ F(0) = \{ \emptyset, \Omega  \}. $ The variable

$$
\theta(n) :=  \left( \sum_{k=1}^n {\bf E} \xi_k^2/F(k-1) \right)^{1/2}, \eqno(2.7)
$$
so that $ S_1(p) = |\theta(n)|_p,  $ is named in \cite{Osekowski1} by "conditional square function"
of our martingale and the variable $ \theta^2(n) \ - $ by
"quadratic (predictable) characteristic" in the review \cite{Peshkir1}. \par
 We deduce using Iensen and triangle inequalities  taking into account the restriction $ p \ge 4: $

$$
\theta^2(n) = \sum_{k=1}^n {\bf E} \xi_k^2 /F(k-1),
$$

$$
|\theta^2(n)|_{p/2} \le  \sum_{k=1}^n | \ {\bf E} \xi_k^2 /F(k-1) \ |_{p/2} \le \sum_k |\xi_k|_p^2 = \sum_{k=1}^n \mu_k^2(p). \eqno(2.8)
$$
 Since

$$
|\theta^2(n)|_{p/2} = |\theta(n)|_p^2,
$$
we ascertain

$$
 S_1(p) = |\theta(n)|_p \le \sqrt{ \sum_{k=1}^n \mu_k^2(p) }. \eqno(2.9)
$$

\vspace{3mm}

 Let us estimate now the value $ S_2(p). $ This evaluate is simple:

$$
S_2^p(p) = {\bf E} \left( \sum_k |\xi_k|^p  \right) = \sum_k |\xi_k|_p^p = \sum_{k=1}^n \mu_k^p(p),\eqno(2.10)
$$

$$
S_2(p) \le \left(  \sum_{k=1}^n \mu_k^p(p) \right)^{1/p} \le  \sqrt{ \sum_{k=1}^n \mu_k^2(p) }.  \eqno(2.11)
$$

  Thus,

$$
 \left| \ \sum_{k=1}^n \xi_k \ \right|_p \le K_{Os} \cdot \frac{p}{\ln p}  \cdot \sqrt{ \sum_{k=1}^n \mu_k^2(p) }. \eqno(2.12)
$$

\vspace{3mm}

{\bf 2.}  Since the sequence $ \{  b(i) \} $ is non-random, the random sequence $ \{\xi^{(b)}(k) \} := \{ b(k) \cdot \xi(k)  \} $ is also
a sequence of martingale differences relative at the same filtration. We apply the last inequality (2.11) for the martingale differences
$ \{\xi^{(b)}(k) \}:  $

$$
 \left| \ \sum_{k=1}^n b(k) \xi_k \ \right|_p \le K_{Os} \cdot \frac{p}{\ln p}  \cdot \sqrt{ \sum_{k=1}^n  b^2(k)  \mu_k^2(p) }, \eqno(2.13)
$$
and we obtain after taking supremum over $  \vec{b} \in B: $

$$
U(p;1,n) = \sup_{b \in B} | \ Q(1,n, b, \{ \xi(\cdot) \}) \ |_p \le  K_{Os} \cdot \frac{p}{\ln p} \cdot \sup_k \mu_k(p),
$$
or equally

$$
U(p;1,n)  \le \gamma(1) \cdot \frac{p}{\ln p } \cdot V_1(p). \eqno(2.14)
$$

\vspace{3mm}

{\bf 3. Remark 2.1.}  We deduce  as a particular case choosing in (2.13) $ b(k) = 1/\sqrt{n}: $

$$
n^{-1/2} \left| \ \sum_{k=1}^n  \xi_k \ \right|_p \le K_{Os} \cdot \frac{p}{\ln p}  \cdot \sqrt{ n^{-1} \sum_{k=1}^n  \mu_k^2(p) }, \eqno(2.15)
$$
which is some generalization of the classical Rosenthal's inequality on the martingale case and in turn is a slight
simplification of the A.Osekovski result. \par

 In turn,

$$
\sup_n \left[ n^{-1/2} \left| \ \sum_{k=1}^n  \xi_k \ \right|_p \right] \le K_{Os} \cdot \frac{p}{\ln p}  \cdot  \sup_{k} \mu_k(p). \eqno(2.15b)
$$

\vspace{3mm}

{\bf 4. Induction step }  $ d \to d+1 $  is completely analogous to one in the article \cite{Ostrovsky7}, section 3, and may be omitted.\par

\vspace{3mm}

{\bf 5. Remark 2.2; an example.} Suppose $ d = \const \ge 2, \ n \ge d + 1. $
  We deduce  as a particular case choosing

$$
b(\vec{i}) = 1/\sqrt{n(n-1) \ldots (n-d+1)} \sim n^{-d/2}: \hspace{6mm} n^{-d/2} |R(d)|_p =
$$

$$
 n^{-d/2}| Q(d,n, \vec{1})|_p = n^{-d/2} \left| \ \sum_{\vec{i} \in I(d,n) }  \xi(\vec{i}) \ \right|_p \le
C(d) \cdot \frac{p^d}{\ln^d p} \cdot \left[  \sup_{i,m} \mu_{i,m}(p) \right]^d. \eqno(2.16)
$$

 \vspace{4mm}

\section{Independent case. }

\vspace{4mm}

The reasoning is basically at the same as in the last section. We will use the famous Rosenthal's inequality
(more exactly, a consequence of this inequality) \cite{Rosenthal1} instead the Osekowski's estimate:

$$
n^{-1/2} \left|\sum_{k=1}^n  \xi_k \right|_p \le K_{R} \cdot \frac{p}{\ln p}  \cdot \sqrt{ n^{-1} \sum_{k=1}^n  \mu_k^2(p) }, \eqno(3.1)
$$
where now $ \{\xi_k \} $ is the sequence of the centered independent random variables with finite $ p^{th} $ moment,
$  K_R $ is the Rosenthal's constant. This estimate is non-improvable. \par
 The exact value of this constant is obtained in \cite{Ostrovsky9}:
 $$
 K_R \approx 1.77638/e \approx 0.6535.
 $$

\vspace{3mm}

 Define the announced sequence $ \kappa = \kappa(d), \ d = 1,2,\ldots $  as follows: $ \kappa(1) := K_R $ and by the following recursion

$$
\kappa(d+1) = \kappa(d) \cdot K_{Os} \cdot \left( 1 + \frac{1}{d} \right)^d. \eqno(3.2)
$$
 Since

$$
\left( 1 + \frac{1}{d} \right)^d \le e,
$$
we conclude

$$
\kappa(d) \le K_R \cdot (K_{Os} \cdot e)^{d-1}, \ d = 1,2,\ldots. \eqno(3.3)
$$

\vspace{3mm}

{\bf Theorem 3.1.} {\it  Let the sequence   } $ \kappa(d) $ {\it  be defined in (3.3). Then in the considered here independent case }

$$
U(p;d,n) \le \kappa(d) \cdot \frac{p^d}{(\ln p)^d } \cdot W_d(p) =
\kappa(d) \cdot \frac{p^d}{(\ln p)^d } \cdot \prod_{m=1}^d \mu_m(p). \eqno(3.4)
$$

\vspace{3mm}

{\bf Remark 3.1.} Let us emphasise the difference between martingale and independent cases. This difference  is except  the coefficient
but  in the factors $  V_d(p) = \prod_{m=1}^d \mu_m(d \cdot p) $  and   $ W_d(p) = \prod_{m=1}^d \mu_m(p). $ \par
 It is clear that there are many examples  when $ W_d(p) < \infty  $   but $ V_d(p) = \infty. $ \par

\vspace{4mm}

\section{Exponential bounds for tails of polynomial martingales. }

\vspace{4mm}

 We intend in in this section to obtain the {\it exponential } bounds for tails of distribution for  the r.v. $ Q(d,n) $
through its (obtained) moments estimates.  We can consider only the martingale case (section 2). \par

\vspace{3mm}

{\bf Theorem 4.1.}  Suppose that the described below sequence of the mean zero martingale differences $ \{ \xi(i,m)  \}  $ satisfies
the restriction

$$
\sup_{i,m} \max( {\bf P}( \xi(i,m) \ge x  ), {\bf P}( \xi(i,m) \le - x  )   ) \le \exp \left( - C_1 x^{q} \ (\ln x)^{- q \ r} \right),\eqno(4.1)
$$

$$
x > e, \ C_1 = \const > 0, \ q = \const > 0, \ r = \const. \eqno(4.1a)
$$

 Then

$$
\sup_{b \in B} \max( {\bf P}( Q(d,n,b) \ge x  ), {\bf P}( Q(d,n,b) \le - x)) \le
$$

$$
\exp \left\{ - C_2 \ x^{q/(dq + 1) } \ (\ln x)^{- q(r-d)/( dq+ 1 ) } \right\},  \ x > e.  \eqno(4.2)
$$

 \vspace{3mm}

 {\bf Proof.} It follows from the theory of the so-called Grand Lebesgue spaces \cite{Kozatchenko1},
\cite{Ostrovsky3}, chapter 1, section 1.8,  \cite{Ostrovsky7}  that the inequality (4.1) is equivalent to the
finiteness of the following norm

$$
\sup_{i,m} \sup_{p \ge 4}  \left[ |\xi(i,m)|_p \cdot p^{-1/q} \cdot \log^{-r} p   \right] = C_3 < \infty, \eqno(4.3)
$$
or equally

$$
\sup_{i,m}   |\xi(i,m)|_p  \le C_3 \ p^{1/q} \ \log^r p.  \eqno(4.3a)
$$

 We apply the theorem 2.1:

$$
\sup_{b \in B}| Q(d,n,b)|_p \le C_4 \ p^{d + 1/q} \ [\log p]^{r - d},
$$
which is in turn equivalent to the proposition (4.2). \par

 Note that other exponential bounds for tail of distribution for the r.v. $ Q(d,n,b)  $  under some additional conditions
 is obtained in \cite{Ostrovsky4}. \par

 \vspace{4mm}

\section{Concluding remarks. }

\vspace{4mm}

{\bf A. Examples of lower estimates.}\\

\vspace{3mm}

 Denote in the independent case

$$
K_I(p; d) = \sup_n \sup_{b \in B} \sup_{{\xi(i,m): |\xi(m,i)|_p < \infty}} \left[ \frac{Q(d,p)}{\prod_{m=1}^d \mu_m(p)} \right].\eqno(5.1)
$$
where  the last $ "\sup" $ is calculated over all the sequences of the centered {\it independent} variables  $ \{  \xi(i,m) \} $
satisfying the condition $ |\xi(i,m)|_p < \infty. $    We obtained

$$
K_I(p; d) \le \frac{C_0(d) \ p^d}{ \ln^d p}, \ C_0(d) = \const > 0.
$$

Our new statement:

$$
K_I(p; d) \ge \frac{C(d) \ p^d}{ \ln^d p }, \ C(d) = \const > 0. \eqno(5.2)
$$

 Proof is very simple. The moment estimations are derived in \cite{Kallenberg1} for
the symmetrical polynomials on mean zero independent identical symmetrically
distributed variables, i.e. particular case for us, for which it is proved that

$$
\frac{|Q(d,n)|_p}{ \sqrt{\Var Q(d,n)} \ \mu^d(p)} \ge  \frac{C_1(d) \ p^d}{ \ln^d p}, \ C_1(d) = \const > 0. \eqno(5.3)
$$

\vspace{4mm}

 Another a more simple example. Let $ n  = 1 $ and  a r.v. $ \eta $ has a Poisson distribution with unit parameter:

$$
{\bf P} (\eta = k) = e^{-1}/k!,
$$
and define $ \xi = \eta - 1, $ then  the r.v. $ \xi $ is centered and

 $$
 p \to \infty \Rightarrow   |\xi|_p \sim p/( e \cdot \ln p).
 $$

 Let also $  \xi_j, \ j = 1,2,\ldots  $ be independent copies of $ \xi. $ Then

$$
\left|\prod_{j=1}^d  \xi_j  \right|_p \sim e^{-d} \frac{p^d}{(\ln p)^d}, \ p \to \infty. \eqno(5.4)
$$

\vspace{3mm}

{\bf B. Estimations for normed variables.} \\

\vspace{3mm}

 Denote

$$
\Psi(b) = \Psi(d,n,b) = \Var( Q(d,n,b)), \ b \in B, \eqno(5.5)
$$
and impose the following condition on the martingale distribution

$$
0 < C_1(d)  \le \sup_n \sup_{b \ne 0} \left[ \frac{\Psi(d,n,b)}{||b||^2} \right] \le C_2(d) < \infty. \eqno(5.6)
$$
 This condition was introduced and investigated in \cite{Ostrovsky4}. \par
We define also a so - called  {\it relative moments} for the r.v. $ \{  \xi(i,m) \} $ under the {\it natural norming:}

$$
\tilde{\mu}_m(p) \stackrel{def}{=} \sup_i \left| \ \xi(i,m)/\sqrt{\Var(\xi(i,m))} \ \right|_p, \eqno(5.7)
$$

 Denote also in the martingale case

 $$
\tilde{V}(p) = \tilde{V}_d(p) \stackrel{def}{=} \prod_{m=1}^d \tilde{\mu}_m(d \cdot p). \eqno(5.8a)
 $$
and  for the independent variables $ \{  \xi(i,m) \} $

 $$
 \tilde{W}(p) = \tilde{W}_d(p) \stackrel{def}{=} \prod_{m=1}^d \tilde{\mu}_m(p). \eqno(5.8b)
 $$

{\bf We will derive as before the  moment estimations of a form }

$$
\tilde{U}(p;d,n) = \tilde{U}(p) \stackrel{def}{=}
$$

$$
\sup_{b \in B} \left| \ Q(d,n, b, \{ \xi(\cdot) \}) /\sqrt{ \Var ( Q(d,n, b, \{ \xi(\cdot) \})  )  } \ \right|_p \le
$$

$$
\tilde{Q}(p;d,n) = \tilde{Q}(p).\eqno(5.9)
$$

 Namely,

$$
\tilde{U}(p;d,n)  \le \gamma(d) \cdot \frac{p^d}{(\ln p)^d } \cdot \tilde{V}_d(p), \eqno(5.10a)
$$
in general (martingale) case  and

$$
\tilde{U}(p;d,n) \le \kappa(d) \cdot \frac{p^d}{(\ln p)^d } \cdot \tilde{W}_d(p), \eqno(5.10b)
$$
for the independent variables; it will be presumed of course the finiteness of the variables $ \tilde{V}(p) $ and $  \tilde{W}(p). $ \par

\vspace{4mm}

 Authors hope that the  last two estimates are more convenient for the practical using. \par

\newpage

{\bf C. Possible generalizations.} \\

\vspace{3mm}

 It is interest by our opinion to generalize our estimates on the {\it predictable} sequence $ b(\vec{i}). $
A preliminary (one-dimensional) result in this direction  see in the article \cite{Ostrovsky0};
see also \cite{Choi1}.\par


\vspace{4mm}

\begin{thebibliography}{99}

\bibitem{Burkholder1}
{\sc D. L. Burkholder.} {\it Distribution Functions Inequalities for Martingales.} Ann. Probab., {\bf 1}, (1973),
19-42.

\bibitem{Burkholder2}
{\sc D. L. Burkholder.} {\it Explorations in Martingale Theory and its Applications.} Ecole d'Ete de
Probabilities de Saint-Flour XIX, 1989, pp. 1-66, Lecture Notes in Math., 1464, Springer, Berlin, 1991.

\bibitem{Choi1}
{\sc Choi K.P.} {\it Some sharp inequalities for martingale transform.} Trans. Amer.
Math. Soc., 307, No 1, (1988), 279-300, MR0936817.

\bibitem{Fan1}
{\sc Fan X., Grama I, and Liu Q.}
{\it Large deviations for martingales with exponential condition.}
arXiv:1111.1407 [math.PR] 6 Nov 2011.

\bibitem{Grama1}
{\sc Grama I.G.}  {\it On moderate deviations for martingales.}
 Ann. Probab. Volume 25, Number 1, (1997), 152-183.

\bibitem{Grama2}
{\sc Grama I.G. and Haeusler E.} (2000.) {\it Large deviations for martingales via Cramer's
method.} Stochastic Process. Appl. 85, 279–293.

\bibitem{Hitczenko1}
{\sc P. Hitczenko.} {\it Best constants in martingale version of Rosenthal inequality.} Ann. Probab. 18
No. 4 (1990), 1656-1668.

\bibitem{Kallenberg1}
{\sc Kallenberg Olav and Sztencel Rafal.} {\it Some dimension-free features
of vector-valued martingales.} Probability Theory Related Fields.
1991, {\bf 88,}  215-247.

\bibitem{Kozatchenko1}
{\sc Kozatchenko Yu. V., Ostrovsky E.I.} {\it The Banach Spaces of
random Variables of subgaussian type.} Theory of Probab. and Math.
Stat. (in Russian). Kiev, KSU, 32, 43-57, (1985).

\bibitem{Lesign1}
{\sc Lesign E., Volny D.} {\it Large deviations for martingales.}
Stochastic Processes and their Applications.  {\bf 96,} 143-159. \  (2001).

\bibitem{Li1}
{\sc Li Y.} ( 2003.) {\it  A martingale inequality and large deviations.} Statist. Probab. Lett., {\bf 62,}
317-321.

\bibitem{Liu1}
{\sc Liu Q., Watbled F.} {\it Exponential inequalities for martingales and asymptotic properties of the
free energy of directed polymers in random environment.}
arXiv:0812.1719v1 [math.PR] 9 Dec 2008.

\bibitem{Osekowski1}
{\sc Osekowski A.} {\it A Note on Burkholder-Rosenthal Inequality.  }
Bull. Polish Academy of Science, Math., {\bf 60}, (2012), 177-185.

\bibitem{Ostrovsky0}
{\sc Ostrovsky E., Sirota L.} {\it  Moment and tail estimates for martingales and martingale transform,
with application to the martingale limit theorem in Banach spaces.}
arXiv:1206.4964v1 [math.PR] 21 Jun 2012

\bibitem{Ostrovsky3}
{\sc Ostrovsky E.I.} {\it Exponential estimations for Random Fields and
     its applications.} (in Russian). 1999, Moskow-Obninsk, Russia, OINPE.

\bibitem{Ostrovsky4}
{\sc  Ostrovsky E. } {\it Bide-side exponential and moment inequalities  for
      tails of distribution of Polynomial Martingales.}
      arXiv: math.PR/0406532 v.1  Jun. 2004

\bibitem{Ostrovsky7}
{\sc Ostrovsky E. and  Sirota L. } {\it Moment and tail inequalities for polynomial martingales.
 The case of heavy tails.}
arXiv:1112.2768v1 [math.PR] 13 Dez 2011

\bibitem{Ostrovsky8}
{\sc Ostrovsky E. and  Sirota L.} {\it Non-improved uniform   tail estimates for normed sums of
independent random variables with heavy tails, with applications.}
arXiv:1110.4879v1 [math.PR] 21 Oct 2011.

\bibitem{Ostrovsky9}
{\sc Ostrovsky E. and  Sirota L.} {\it Schl\"omilch and Bell series for Bessel's functions, with probabilistic applications.}
arXiv:0804.0089v1 [math.CV] 1 Apr 2008


\bibitem{Peshkir1}
{\sc Peshkir G., Shirjaev A.N.} {\it The Khintchine inequalities and martingale
expanding sphere of their action.} Russian Math. Surveys; {\bf 50,} 5, 849-904, (1995).


\bibitem{Volny1}
{\sc Volny D.} {\it Approximating martingales and the central limit theorem for strictly  stationary
processes.} Stoch. processes ant their applic., {\bf 44}, (1993),  41-74.

\bibitem{Ratchkauskas1}
{\sc Ratchkauskas A.} {\it Large deviations for martingales with some applications.}
Acta Applicandae Mathematicae, Volume 38, Number 1 (1995), 109-129, DOI: 10.1007/BF00992617

\bibitem{Rosenthal1}
{\sc Rosenthal H.P.}  {\it On the Subspaces of $ L_p, \ (p > 2) $ spanned by Sequences of
independent Variables. } Israel J. Math., 1970, V.3, pp. 253-273.

\end{thebibliography}
\end{document}